\documentclass[reqno,10pt]{article}
\usepackage{graphicx}
\usepackage{amsmath}
\usepackage{amssymb}
\usepackage{amsthm}
\usepackage{enumitem}
\usepackage{bbm}
\usepackage[l2tabu,orthodox]{nag}
\usepackage[all,error]{onlyamsmath}
\usepackage[strict]{csquotes}
\usepackage{url}
\usepackage{tikz}
\usepackage{lipsum}

\newtheorem{theorem}{Theorem}[section]
\newtheorem{lemma}[theorem]{Lemma}

\theoremstyle{definition}
\newtheorem{definition}[theorem]{Definition}
\newtheorem{example}[theorem]{Example}
\theoremstyle{remark}
\newtheorem{remark}[theorem]{Remark}

\numberwithin{equation}{section}

\numberwithin{equation}{section}
\newcommand{\be}{\begin{equation}}
\newcommand{\ee}{\end{equation}}
\newcommand{\ba}{\begin{aligned}}
\newcommand{\ea}{\end{aligned}}

\newcommand{\N}{{\mathbb N}}
\newcommand{\Z}{{\mathbb Z}}
\newcommand{\R}{{\mathbb R}}

\def\csi1{\circ\sigma^{-1}}

\newcommand\blfootnote[1]{%
  \begingroup
  \renewcommand\thefootnote{}\footnote{#1}%
  \addtocounter{footnote}{-1}%
  \endgroup
}
\newcommand{\B}{{\mathcal B}}

\newenvironment{OMabstract}{\noindent\textbf{Abstract.} }{\medskip}
\newenvironment{OMsubjclass}{\noindent\textbf{Mathematics Subject Classification (2020):} }{\medskip}
\newenvironment{OMkeywords}{\noindent\textbf{Keywords:}  }{\medskip}
\date{\vspace{-5ex}}

\begin{document}

\author{Shrey Sanadhya} 
\title{Coboundaries of commuting Borel automorphisms}
\maketitle


\begin{OMabstract}
    We show that if $S,T$ are two commuting automorphisms of a standard Borel space such that they generate a free Borel $\Z^2$-action then $S$ and $T$ do not have same sets of real valued bounded coboundaries. We also prove a weaker form of Rokhlin Lemma for Borel $\Z^d$-actions. 
\end{OMabstract}

\begin{OMkeywords}
    coboundries, Rokhlin Lemma, Borel $\mathbb{Z}^d$-action.
\end{OMkeywords}

\begin{OMsubjclass}
    37A40, 37A99, 37B99.
\end{OMsubjclass}


\section{Introduction}\label{sect Intro}

Let $\Gamma$ be a countable Borel automorphism group of a standard Borel space $(X,\B)$. A real valued Borel map $\alpha : \Gamma \times X \to \mathbb{R}$ is called a \textit{cocycle} of $\Gamma$ if it satisfies the following cocycle identity
\begin{equation}\label{cocy}
    \alpha(\gamma_1\gamma_2, x) = \alpha(\gamma_1, \gamma_2 x) 
+ \alpha(\gamma_2, x),
\end{equation} for all $\gamma_1,\gamma_2 \in \Gamma$ and $x \in X$. A cocycle $\alpha(\gamma, x)$ is called a \textit{coboundary} if there 
exists a Borel function $c : X \to \mathbb{R}$ such that $\alpha(\gamma, x) = f(\gamma x) - f(x)$. \blfootnote{This article is a part of the author’s Ph.D. thesis, written under the supervision of Sergii Bezuglyi and Palle E. T. Jørgensen, University of Iowa.}

Cocycles have played a central role in ergodic theory, where they are defined upto zero measure sets : if $\Gamma$ is a countable discrete group acting on a standard probability space $(X,\mu,\B)$ by non-singular transformations, then relation \eqref{cocy} must hold $\mu$-a.e. The study of cocycles have given rise to constructions (for e.g. skew products) which have been crucial in the classification of dynamical systems up-to orbit equivalence. Moreover, in context of ergodic theory, the study of cocycles has resulted in the classification which is even finer than orbit equivalence (see for e.g., \cite{BezuglyiGolodets_1991}, \cite{GolSinelsh1994}, \cite{GolodetsSinelshchikov_1987},\cite{Hamachi_2000}, \cite{FeldmanSutherlandZimmer1989}). Apart from this circle of problems, cocycles proved their usefulness in representation theory, theory of groupoids, operator algebras etc. A complete list of papers devoted to the study of cocycles of dynamical systems, containing all crucial contributions to the field, is too long to mention. We list several key reference which include the seminal work of Ramsay \cite{Ramsay_1971}, Moore \cite{Moore_1970}, Feldman \cite{FeldmanMoore_1977}, Zimmer \cite{Zimmer_1984},
Schmidt \cite{Schmidt_1977, Schmidt_1990}.

This paper is focused on the study of Borel cocycles (in particular the structure of set of coboundaries) of Borel dynamical systems. There exist many similarities between ergodic theory and Borel dynamics, but the fact that, in Borel dynamics, there is no prescribed measure on the underlying space makes these two theories substantially different. Following the work in ergodic theory, it is natural to study problem involving Borel cocycles of dynamical systems. For example, it would be important to know if a version of the classical notion of ratio sets (see \cite{Schmidt_1977} for details) can be studied for Borel cocycles. Similarly, it would be interesting to study a version of Mackey range (see \cite{Mackey_1966}) in context of Borel cocycles. These ideas have proved to be extremely useful in classification of automorphism groups in ergodic theory. We should also mention that the properties of Borel cocycles have been studied in the following papers : \cite{Becker_2013}, \cite{ConzeRaugi_2009}, \cite{Danilenko1998}, \cite{FeldmanMoore_1977}, \cite{Miller_2006}, \cite{Miller_2008}, \cite{Conley_Miller_2017}, \cite{Miller_2020} and more.

A key motivation to study Borel cocycles stems from the theory of orbit equivalence of Borel dynamical systems. The study of orbit equivalence in Borel dynamics is equivalent to the study of isomorphism of the corresponding equivalence relations (generated by orbits) in the descriptive set theory. Thus, the study of Borel cocycles is closely connected to the notion of countable Borel equivalence relation (CBER) which has been extensively studied in the descriptive set theory. We encourage interested readers to 
\cite{BeckerKechris_1996}, \cite{DoughertyJacksonKechris_1994},  
\cite{JacksonKechrisLouveau_2002}, \cite{Hjorth_2000}, 
\cite{Kechris_1995}, \cite{KechrisMiller_2004}, \cite{Nadkarni_2013},
\cite{Varadarajan_1963}, and 
\cite{Weiss_1984}, where one can find connections between the theory of orbit equivalence in Borel dynamics and the descriptive set theory.

We fix the main setting for the paper: Let $T$ be an aperiodic Borel automorphism of a standard Borel space $(X,\mathcal{B})$. Every real valued Borel function $f : X \to \R$ defines a cocycle 
$a: \mathbb{Z} \times X \rightarrow \R$ (of the Borel $\Z$-action corresponding to $T$) by the formula   
\begin{equation}\label{eq z-cocycle}
    a(j,x) = \begin{cases}
f(x)+ f(Tx)+...+f(T^{j-1}x), & \quad j\geq 1 \\ 
0,  & \quad j = 0\\ 
-f(T^{-1}x) - f(T^{-2}x) - ... -f(T^{j}x), & \quad j\leq-1,
\end{cases}    
\end{equation}

Conversely, if  $a : \mathbb Z \times X \to \R$ is a cocycle of the group 
$\{T^n, n \in \mathbb{Z}\}$, then it is completely determined by the
function $f(x) = a (1, x)$. Moreover, the properties of the cocycle $a(j, x)$
are represented in terms of the function $f$.

A real valued Borel map $f$ on $X$ is called a (Borel) \textit{coboundary} for $T$ if there exists a Borel function $g$ such that
\begin{equation}\label{eq_cob}
    f(x) = g(x) - g(Tx)
\end{equation} for every $x \in X$. Two Borel functions $f$ and $h$ are called \textit{cohomologous} if $f-h$ is a coboundary. We denote by $Cob(T)$ the set of real valued bounded (Borel) coboundaries of $T$. 

This paper is inspired by the following theorem of Kornfeld \cite{Kornfeld_1997} from ergodic theory :

\begin{theorem} [Kornfeld \cite{Kornfeld_1997}]\label{Ko} Suppose $\sigma$ and $\tau$ are two commuting invertible ergodic measure preserving transformations of a non-atomic probability space $(X,\B,\mu)$. They have the same coboundaries if and only if $\sigma = \tau^{\pm1}$. 
\end{theorem}The coboundries in Theorem $\ref{Ko}$ are measurable real valued functions that satisfy $(\ref{eq_cob})$ for $\mu$ a.e. $x \in X$. A similar result in the context of Cantor dynamics is due to N. Ormes (Theorem $\ref{Or}$). We refer the reader to \cite{Ormes_2000} for details. Coboundaries considered in Theorem $\ref{Or}$ are continuous.

\begin{theorem} [Ormes \cite{Ormes_2000}\label{Or}] Let $(X,T)$ and $(Y,S)$ be two Cantor minimal systems. There is an orbit equivalence $h : X \rightarrow Y$ which induces a bijection from the set of real $S$-coboundaries to the set of real $T$-coboundaries if and only if $S$ and $T$ are flip conjugate, i.e. $S$ is conjugate to $T$ or $S$ is conjugate to $T^{-1}$.

\end{theorem} Following Kornfeld and Ormes, the key question we ask in this note is as follows: To what extent does the set $Cob(T)$ determine $T$? We adopt Kornfeld's ideas to the domain of Borel dynamics and prove the following theorem :

\begin{theorem}\label{thm main 1} Let $S,T$ be two commuting automorpisms of standard Borel space such that they generate a free Borel $\Z^2$-action. Then $S$ and $T$ do not have same sets of real valued bounded coboundaries.
\end{theorem} 

\begin{remark}
Note that $S = T^{-1}$ implies $Cob(T) = Cob(S)$. To see this, let $f \in Cob(T)$, thus $f(x) = h(x)-h(Tx)$ for some Borel function $h$ and every $x \in X$. Then $f$ can also be written as $f(x) = h_1(x) - h_1(Sx)$ where $h_1 = -h \circ T$. Theorem $\ref{thm main 1}$ gives us a criterion to determine when two commuting transformations do not have same set of bounded coboundries. 
\end{remark}

A key component in the proof of Theorem $\ref{thm main 1}$ is a \textit{two-dimensional version} of the following theorem (also see \cite[Lemma 3.3]{BezuglyiDooleyKwiatkowski_2006}) proved for a single automorphism:

\begin{theorem} [Weiss \cite{Weiss_1984}] \label{Marker_0} Let $S$ be an aperiodic Borel automorpisms on a standard Borel space $(X, \B)$. Then there exists a decreasing sequence $A_n$, $n \in \N$, of Borel sets such that

\noindent $(i)$ for each $n$, $\underset{i= - \infty}{\overset{\infty}{\bigcup}} S^i A_n = \underset{i= - \infty}{\overset{\infty}{\bigcup}} S^i (X - A_n)$ = X,

\noindent $(ii)$ for each $n$, the sets $A_n, S(A_n), S^2(A_n), ...S^{n-1} A_n$ are pairwise disjoint,

\noindent $(iii)$ the intersection $\underset{n=1}{\overset{\infty}{\bigcap}} A_n  = \emptyset$. 

\end{theorem} The nested sequence of Borel sets $A_n$, satisfying conditions $(i)-(iii)$ of Theorem $\ref{Marker_0}$, is called a \textit{vanishing sequence of markers}. We prove a finite-dimensional version of this theorem (Theorem $\ref{d-Marker}$). 

Let $\Gamma$  denote a Borel action of a countable group (see Definition $\ref{Borel-act}$), then a Borel set $A$ is a \textit{complete section} for $\Gamma$ if every $\Gamma$-orbit intersects the set $A$. Following theorem is a finite-dimensional version of Theorem $\ref{Marker_0}$. 

\begin{theorem}\label{d-Marker} Let $\Gamma$ be a free Borel $\Z^d$-action generated by commuting Borel automorphisms, $\{T_i\}$, $i=1,2,..,d$, of a standard Borel space $(X,\B)$. Then for any $d$-tuple $(t_1,...,t_d)$ of positive integers, there exists a Borel set $A$, such that 

\vspace{2mm}

\noindent $(i)$  $A$ is a complete section for $\Gamma$,

\vspace{2mm}

\noindent $(ii)$ the sets $T_1^{k_1}...T_d^{k_d} A$, $0\leq k_i < t_i$, $i=1,2,..,d$, are pairwise disjoint.

\end{theorem} Theorem $\ref{d-Marker}$ can also be considered as a weaker version of \textit{Rokhlin Lemma} for Borel $\Z^d$-action. We should mention that Theorem $\ref{d-Marker}$ can be derived from the work of S. Gao and S. Jackson (see Theorem 3.1 in \cite{GaoSuJackson_2015}) where the authors used descriptive set theoretic methods (big-marker-little-marker method) to prove the theorem. In this paper, we provide a different and purely dynamical proof for Theorem $\ref{d-Marker}$. 

The \textit{outline} of the paper is as follows. In Section $\ref{sect Prilim}$, we provide basic definitions and preliminary results about groups of Borel automorphisms. In Section $\ref{Sec weak}$, we provide a proof of Theorem $\ref{d-Marker}$. In Section $\ref{cob}$, we use Theorem $\ref{d-Marker}$ (for $d=2$) to prove Theorem $\ref{thm main 1}$.

\textit{Notation and Terminology:} Here are a few remarks about the notation and terminology of our results in this paper. We have preferred to use the terminology which is traditional for ergodic theory and dynamical systems. Hence, our main objects of study are Borel automorphism groups and not the equivalance relations. But when it is convenient, we have also used the language of countable Borel equivalence relations (CBER). For the ease of readers, we provide basic definitions and terminology from descriptive set theory in Section $\ref{sect Prilim}$. Throughout the paper, we fix the following \textbf{notation}:

\begin{itemize}

\item $(X,\mathcal{B})$ is a standard Borel space with the 
$\sigma$-algebra of Borel sets $\mathcal{B}= \mathcal{B}(X)$.

\item A one-to-one Borel map $T$ of space $(X,\mathcal{B})$ onto itself 
is called a Borel automorphism of $X$. In this paper the term 
"automorphism"  means a Borel automorphism of $(X,\mathcal{B})$.

\item In this paper $Cob(T)$ denotes the set of real valued bounded coboundries (see $(\ref{eq_cob})$) of  Borel automorphism $T$. 

\item $Aut(X,\mathcal{B})$ is the group of all Borel automorphisms of 
$X$ with the identity map ${\mathbb I} \in Aut(X, \mathcal{B})$.

\end{itemize}

\section{Preliminaries}\label{sect Prilim} In this section we provide the basic definitions from Borel dynamics and descriptive 
set  theory.

\noindent \textbf{Automorphisms of standard Borel space.} Let $X$ be a separable completely metrizable topological space (also called \textit{Polish space}). Let $\B$ be the $\sigma$-algebra generated by the open sets in $X$. Then we call the pair $(X,\B)$ a \textit{standard Borel space}. Any two uncountable standard Borel spaces
are isomorphic. Any countable subgroup $\Gamma$ of $Aut(X, \mathcal B)$ is called a \textit{Borel automorphism group}, and the pair $(X,\Gamma)$ is referred to as a \textit{Borel dynamical system}.  

\begin{definition}\label{Borel-act} Let $G$ be a countable group with  identity $e$. A \textit{Borel action} of the group $G$ on $(X, \mathcal B)$ is a group homomorphism $\rho : g \rightarrow \rho_g : G
 \rightarrow Aut(X, \mathcal{B})$. In other words, for each $g\in G$, $\rho_g : X \rightarrow X$ is a Borel automorphism such that 
(i) $\rho_{gh}(x)=\rho_g(\rho_h (x))$ for every $h \in G$ and  (ii) $\rho_{e}(x)=x$ for every $x \in X$.

\end{definition} Note that $\rho(G)$ is a countable subgroup of $Aut (X,\mathcal{B})$. If, for some $x \in X$, the relation  $\rho_g(x) = x$  implies $g = e$, then $\rho$ is called a \textit{free action} of $G$. In this case, the group homomorphism $\rho$ is injective. Any Borel automorphism $T \in Aut(X,\B)$ gives rise to a Borel action of group $\mathbb{Z}$ (also referred as a Borel $\Z$-action) by identifying $k \in \mathbb{Z}$ with $T^k \in Aut(X,\B)$. A Borel automorphism $T$ is called \textit{periodic at a point} $x \in X$ if there exists $k \in \mathbb N$ such that $T^kx = x$. The least such $k$ is called the \textit{period} of $T$ at $x$. A Borel automorphism $T$ is called \textit{aperiodic} if every $T$-orbit is countably infinite. In this paper we work with free Borel $\Z^d$-actions generated by $d$ commuting automorphisms. We make this notion precise in the following definition. For simplicity we use $d=2$.  
 
\begin{definition}\label{gen} Let $S$ and $T$ be two commuting Borel automorphism we say that they generate a free Borel action $\rho(\Z^2)$ of group $\Z^2$ where $\rho_{(n,m)} : = S^n T^m = T^m S^n \in Aut(X,\B)$, if for every $(n,m) \in \Z^2$ we have  $\rho_{(n,m)} (x) \neq x$ for all $x \in X$ and any pair of integers $(m,n) \neq (0,0)$.

\end{definition} Note that the above definition immediately implies that $S$ is not a power of $T$ for any $x \in X$. To see this, assume there exist $x \in X$ such that $S x= T^p x$ for some $p\in \Z$. Hence $S T^{-p} x = x$, which contradicts the assumption that $\rho$ is a free action.

\textit{Countable Borel equivalence relation (CBER)}:  An equivalence relation $E$ on $(X, \mathcal{B})$ is called Borel if it is a Borel subset of the product
space $E \subset X \times X$, where $X \times X$ is equipped with the Borel $\sigma$-algebra $\mathcal{B} \times \mathcal{B}$. It is called countable if every equivalence class 
$[x]_E := \{y \in X : (x,y) \in E\}$ is countable for all $x \in X$. If $B$ is a Borel set, then $[B]_E$ denotes the saturation of $B$ with respect to the equivalence relation $E$, i.e., $[B]_E$ contains the entire class $[x]_E$ for every $x\in B$. In the study of Borel dynamical systems the theory of CBER plays an important role as it provides a link between descriptive set theory and Borel actions (see Theorem $\ref{FM}$).

Let $\Gamma$ be a Borel automorphism group of $(X,\B)$, then the \textit{orbit equivalence relation} generated by the action of $\Gamma$ on $X$ is given 
\begin{equation*}
    E_X(\Gamma) = \{(x,y) \in X \times X: x=\gamma y
    \  \mbox{for\ some}\ \gamma \in \Gamma\}.
\end{equation*} Note that $E_X(\Gamma)$ is a CBER. We will call an equivalence relation $E$ \textit{periodic at a point} $x \in X$ if the equivalence class $[x]_E$ is finite. Similarly, an equivalence relation $E$ is \textit{aperiodic at} $x \in X$ if the equivalence class $[x]_E$ is countably infinite. If every $E$-class is countably infinite we will say that the equivalence relation $E$ is \textit{aperiodic}. In this paper we work with aperiodic CBERs. The study of aperiodic CBERs in itself is an area of immense importance in descriptive set theory. We refer our readers to \cite{Kechris_2019} for an up-to-date survey of the theory of countable Borel equivalence relations.
    
The following theorem shows that all CBERs come from Borel actions of countable groups.

\begin{theorem} [Feldman--Moore \cite{FeldmanMoore_1977}] \label{FM}  Let $E$ be a
 countable Borel equivalence relation on a standard Borel space $(X,
 \mathcal{B})$. 
Then there is a countable group $\Gamma$ of Borel automorphisms of 
$(X,\mathcal{B})$ such that $E = E_X(\Gamma)$.
\end{theorem} 

\begin{definition}
A Borel set $C$ is called a \textit{complete section} for an equivalence relation $E$ on $(X,\B)$ if every $E$-class intersects $C$, in other words $[C]_E = X$. Let $\Gamma \in Aut(X,\B)$ be a countable Borel automorphism group. Then we will denote by $\mathcal{C}_\Gamma$ the collection of Borel subsets $C$ such that $C$ and $X \setminus C$ both are complete section for $E_X(\Gamma)$.

\end{definition} If a complete section intersects each $E$-class exactly once then it is called a \textit{Borel transversal}. An equivalence relation $E$ which admits a Borel transversal is called \textit{smooth.} Equivalently, one can say that an equivalence relation  $E$ on a standard Borel space $(X, \B)$ is \textit{smooth} if there is a Borel function $f: X \rightarrow Y$, where $Y$ 
is a standard Borel space, such that $(x,y) \in E \Longleftrightarrow f(x)= f(y)$. We remark that in this paper we will deal only with \textit{non-smooth} CBERs.

\textit{Full group of automorphisms.} 
For a countable subgroup $\Gamma$ of $ Aut(X, \mathcal{B})$, we denote by $\Gamma x$ the orbit  $\{\gamma x : \gamma \in \Gamma \}$  of $x$ with  respect to $\Gamma$. We say that $\Gamma$ is a \textit{free} group of automorphisms  if $\gamma x \neq x$ for every $\gamma \neq e$ and $ x \in X$. The set 
\begin{equation*}
     [\Gamma] = \{R \in Aut (X,\mathcal{B}) : Rx \in \Gamma x, \ \forall x\in X \}
\end{equation*} is called the \textit{full group of automorphisms} generated by $\Gamma$.  The full group generated  by a single automorphism $T \in Aut (X, \mathcal{B})$ is denoted by $[T]$. 

A countable subgroup $\Gamma$ of $Aut (X, \mathcal{B})$ is called 
\textit{hyperfinite} if $\Gamma x = \bigcup _{i=1} ^\infty \Gamma_i x$ for every $x \in X$, where each $\Gamma_i$ is a finite subgroup of $[\Gamma]$ and $\Gamma_i \subset \Gamma_{i+1}$ for all $i$.
Equivalently, a countable Borel equivalence relation $E$ is called 
\textit{hyperfinite} if $E = \bigcup_{n} E_n $ where $E_n \subset E_{n+1}$ for all $n$, where each $E_n$ is a finite Borel sub-equivalence relation of $E$. 

Let $\Gamma_1, \Gamma_2$ be two countable Borel automorphism groups of $(X,\B)$. We say that $\Gamma_1$ and $\Gamma_2$ are \textit{orbit equivalent} (also \textit{o.e.}) if there exists a Borel isomorphism $\phi :X \rightarrow X $ such that $\phi(\Gamma_1 x) = \Gamma_2(\phi(x))$, $\forall x \in X$. In other words $\Gamma_1$ orbit of $x$ is same as the $\Gamma_2$ orbit of $\phi (x)$ for every $x \in X$. The following theorem gives an important characterisation of a hyperfinite CBER. 

\begin{theorem}  [Slaman-Steel \cite{SlamanSteel_1988}, Weiss \cite{Weiss_1984}]\label{thm hyperfinite}
Suppose $E$ is a CBER. The following are equivalent:

$1.$ $E$ is hyperfinite.

$2.$ $E$ is generated by a Borel $\mathbb{Z}$-action.

\end{theorem}

\section{Weak Rokhlin Lemma for Borel $\Z^d$-actions} \label{Sec weak} In this section, we prove Theorem $\ref{d-Marker}$, which can be considered as a weak version of Rokhlin Lemma for free Borel $\Z^d$-action. Before we prove Theorem $\ref{d-Marker}$, we mention its one-dimensional version (Theorem $\ref{Marker1}$) due to B. Weiss (\cite{Weiss_1984}) and a relevant lemma (Lemma \ref{lem decom}). We refer our readers to  \cite[chapter 7]{Nadkarni_2013} for a proof of Theorem $\ref{Marker1}$ and Lemma \ref{lem decom}. 

Let $S$ be an aperiodic Borel automorphism of a standard Borel space $(X, \B)$. A Borel set $W \in \B$ is said to be \textit{wandering} with respect to $S$ if the sets $S^i W$, $i\in \Z$, are pairwise disjoint. We will denote by $\mathcal{W}_S$ (or $\mathcal{W}$ when $S$ is obvious) the sigma ideal generated by all the wandering sets in $\B$.

\begin{definition} A set $A \in \B$ is said to be \textit{decomposable} (mod $\mathcal{W}_S$) if we can
write $A$ as a disjoint union of two Borel sets $C$ and $D$ such that saturation (mod $\mathcal{W}_S$) of $C$, $D$ and $A$ with respect to $S$ are same. In other words,

\begin{equation*}
    \underset{i= - \infty}{\overset{\infty}{\bigcup}} S^i C = \underset{i= - \infty}{\overset{\infty}{\bigcup}} S^i D = \underset{i= - \infty}{\overset{\infty}{\bigcup}} S^i  A \,\,\,\,(\mathrm{mod}\,\, \mathcal{W}_S).
\end{equation*}

\end{definition}

\begin{lemma} \label{lem decom} Let $S$ be a free Borel automorphism on a standard Borel space $(X,\B)$. Then $X$ is decomposable (mod $\mathcal{W}_S$), i.e. there exist a Borel set $A$ such that 

\begin{equation*}
    \underset{i= - \infty}{\overset{\infty}{\bigcup}} S^i A = \underset{i= - \infty}{\overset{\infty}{\bigcup}} S^i (X \setminus A) = X \,\,\,\,(\mathrm{mod}\,\, \mathcal{W}_S). 
\end{equation*} Moreover we can choose the set $A$ such that it does not contain a full $S$ orbit.
\end{lemma}

\noindent\textit{Proof}. See \cite[Lemma 7.24]{Nadkarni_2013} for a proof. \hfill{$\square$}

\begin{theorem} [Weiss \cite{Weiss_1984}] \label{Marker1} Let $S$ be an aperiodic, Borel automorpisms on a standard Borel space $(X, \B)$. Then there exists a decreasing sequence $A_n$, $n \in \N$, of Borel sets such that

\noindent $(i)$ for each $n$, $\underset{i= - \infty}{\overset{\infty}{\bigcup}} S^i A_n = \underset{i= - \infty}{\overset{\infty}{\bigcup}} S^i (X - A_n) = X $, in other words, for each $n$,  $A_n \in \mathcal{C}_{S}$,

\noindent $(ii)$ for each $n$, the sets $A_n, S(A_n), S^2(A_n), ...S^{n-1} A_n$ are pairwise disjoint,

\noindent $(iii)$ the intersection $\underset{n=1}{\overset{\infty}{\bigcap}} A_n =: A_{\infty}$ is a wandering set. 

\end{theorem} 

\noindent\textit{Proof}. See \cite[Theorem 7.25]{Nadkarni_2013} for a proof. \hfill{$\square$}

We will prove Theorem $\ref{d-Marker}$ for $d=2$ (see Theorem $\ref{Marker}$ below). The proof for $d> 2$ follows similarly. In the proof of Theorem $\ref{Marker}$, we will work with \textit{lexicographic order} on a finite index set in $\N_0 \times \N_0$. We define the index set and the ordering in the remark below. 

\begin{remark}\label{Index} For $n,m \in \N$, we denote by $I_{(n,m)}$ the index set containing $nm$ pairs of integers i.e. 
\begin{equation*}
    I_{(n,m)} = \{(0,0),(0,1),...,(n-1,m-1)\}.
\end{equation*} Each element of $I_{(n,m)}$ represents the powers to which $S$ and $T$ are raised. Thus $(i,j) \in I_{(n,m)}$ corresponds to $S^iT^j$. For the rest of the section, we work with the lexicographic order on $I_{(n,m)}$ (denoted by "$\prec$"). For $(i_1,j_1), (i_2,j_2) \in I_{(n,m)}$, we say $(i_1,j_1) \prec (i_2,j_2)$ if,

\noindent $(i)$ $i_1 < i_2$ or

\noindent $(ii)$ $i_1 = i_2$ and $j_1 < j_2.$

\end{remark}

\begin{theorem}\label{Marker} Let $\Gamma$ be a free Borel $\Z^2$-action generated by two commuting aperiodic Borel automorphisms $S,T$. Then for each pair $(n,m) \in \N^2$ there exists a Borel set $A_{(n,m)}$, such that,

\noindent $(i)$ for each pair $(n,m)$, $A_{(n,m)} \in \mathcal{C}_{\Gamma}$,

\noindent $(ii)$ for each $(n,m)$, the sets $S^i T^j A_{(n,m)}$, $(i,j) \in I_{(n,m)}$, are pairwise disjoint.

\end{theorem} 

\noindent\textit{Proof}. A Borel $\Z^d$-action defines a hyperfinite equivalence relation (see \cite{GaoSuJackson_2015}, \cite{Weiss_1984}). Since $\Gamma$ is a $\mathbb{Z}^2$-action, $E_{X}(\Gamma)$ is hyperfinite, hence generated by a Borel $\mathbb{Z}$-action (see Theorem $\ref{thm hyperfinite}$). Thus there exists $R \in Aut(X,\B)$, such that for every $x \in X$, $\Gamma x = \{R^i x : i \in \Z \}$. 

Using lemma $\ref{lem decom}$, choose $A \in \mathcal{C}_{R}$ such that $A$ does not contain a full $R$-orbit. Set $A_{(0,0)} = A$. Let $(n,m)$ be a pair of integers and let $I_{(n,m)}$ be the corresponding index set (see Remark $\ref{Index}$). As mentioned in Remark $\ref{Index}$ above, there are $nm$ elements in the set $I_{(n,m)}$, each corresponding to a power of $S$ and $T$. We consider lexicographic order on the index set. Hence the smallest element of $I_{(n,m)}$ is the pair $(0,0)$ (which corresponds to $S^0T^0$) the next element in the order is the pair $(0,1)$ (which corresponds to $S^0T^1$) and so forth. We want to find $A_{(0,1)} \subset A$ such that, $A_{(0,1)} \in \mathcal{C}_R$ and $A_{(0,1)} \cap T (A_{(0,1)}) = \emptyset$.

Define $\{B_j\}_{j\in \Z \setminus \{0\}}$ a partition of $A = A_{(0,0)}$ as follows,
\begin{equation*}
    B_j = \{x \in A_{(0,0)} : Tx = R^j x\}, \,\,\, j \in \Z \setminus \{0\}. 
\end{equation*} Since $A_{(0,0)}$ does not contain a full $R$-orbit, none of the $B_j$, $j \in \Z \setminus \{0\}$ contains full $R$-orbit. We start with $B_{-1}$. For $x \in B_{-1}$, let $-n(x)$ be the first negative integer $-n$ such that $R^{-n}(x) \notin B_{-1}$. Write
\begin{equation*}
    E_{-k}^{-1} = \{x \in B_{-1} : -n(x) = -k\},\,\,\, k = 1,2,3..
\end{equation*} Since $R^{-n(x)} (x) \notin B_{-1}$, there are two possibilities:

\begin{enumerate}[label={\textup{\arabic*.}}, widest=3, leftmargin=*]

\item$R^{-n(x)} (x) \in X \setminus A_{(0,0)}$ (we call such elements of $B_{-1}$ type $a$ elements) or

\item $R^{-n(x)} (x) \in B_i$, for $i \neq -1$ (we call such elements of $B_{-1}$ type $b$ elements). 

\end{enumerate} We remove type $b$ elements from $B_{-1}$. Thus,
\begin{equation*}
    B_{-1} = B_{-1} \setminus \{ x \in B_{-1} : \, x \,\textrm{is type}\,\, b\},
\end{equation*} and
\begin{equation*}
A_{(0,0)} = A_{(0,0)} \setminus \{ x \in B_{-1}: \, x \,\textrm{is type}\,\, b\}.
\end{equation*} Hence the set $E_{-k}^{-1}$ becomes
\begin{equation*}
E_{-k}^{-1} = \{x \in B_{-1} : -n(x) = -k \,;\,\, R^{-n(x)} (x) \in X \setminus A_{(0,0)} \},\,\,\, k = 1,2,3..
\end{equation*} Put 
\begin{equation} \label{eq 4a}
    A_{-1} = \bigcup_{-k=-1}^{-\infty} R^{(-k+1)} (E_{-k}^{-1}) \subseteq B_{-1}
\end{equation} Note that $A_{-1} \cap R^{-1} (A_{-1})= \emptyset$, and $R^{-1}(A_{-1}) \cap B_i = \emptyset$ for $i \in \Z \setminus \{0\}$, $i \neq -1$.

Now we repeat the same process for $B_1$. For $x \in B_{1}$, let $n(x)$ be the first positive integer $n$ such that $R^{n}(x) \notin B_{1}$. Again there are two possibilities: 

\begin{enumerate}[label={\textup{\arabic*.}}, widest=3, leftmargin=*]

\item $R^{n(x)} (x) \in X \setminus A_{(0,0)}$ (again we denote such elements as type $a$ elements of $B_1$) or
\item $R^{n(x)} (x) \in B_i$, for $i \neq -1,1$ (type $b$ elements of $B_1$).

\end{enumerate} We remove type $b$ elements from $B_{1}$. Thus
\begin{equation*}
B_{1} = B_{1} \setminus \{x \in B_{1}: \, x \,\textrm{is type}\, b\},
\end{equation*} and
\begin{equation*}
A_{(0,0)} = A_{(0,0)} \setminus \{x \in B_{1}: \, x \,\textrm{is type}\, b\}.
\end{equation*} The set $E_{k}^{1}$ is now defined as

\begin{equation*}
E_{k}^{1} = \{x \in B_{1} : n(x) = k\,;\,\, R^{n(x)} (x) \in X \setminus A_{(0,0)} \},\,\,\, k = 1,2,3..
\end{equation*} Put 
\begin{equation} \label{eq 4b}
    A_{1} = \bigcup_{k=1}^{\infty} R^{(k-1)} (E_{k}^{1}) \subseteq B_{1}.
\end{equation} Note that $A_{1} \cap R (A_{1})= \emptyset$, and $R(A_{1}) \cap B_i = \emptyset$ for $i \in \Z \setminus \{0\}$, $i \neq 1$.

We now repeat this process for $B_{-2}$ and then $B_{2}$ and so on. Finally, for every $j \in \Z \setminus \{0\}$ we obtain $A_j \subseteq B_j$ such that, $R^j (A_j) \cap A_j = \emptyset$ and $R^j (A_j) \cap B_i = \emptyset$ where $i \neq j$, $i \in \Z \setminus \{0\}$. Put 
\begin{equation*}
A_{(0,1)} = \underset{i \in \Z \setminus \{0\}}{\bigcup} A_i \subseteq A_{(0,0)},\,\, \textrm{then} \,\, A_{(0,1)} \in \mathcal{C}_R \,\, \textrm{and} \, A_{(0,1)} \cap T (A_{(0,1)}) = \emptyset. 
\end{equation*} Thus, we found a set $A_{(0,1)} \in \mathcal{C}_R$ such that $A_{(0,1)} \cap T (A_{(0,1)}) = \emptyset$. Now we move on to next power in the lexicographic order, i.e. $(0,2)$. Instead of working with $A$ we now work with $A_{(0,1)} \subset A$ and repeat the above procedure to find $A_{(0,2)} \subset A_{(0,1)}$, such that $A_{(0,2)} \in \mathcal{C}_R$ and $A_{(0,2)} \cap T (A_{(0,2)}) = \emptyset$. Thus, we have obtained set $A_{(0,2)} \subset A$ such that $A_{(0,2)} \in \mathcal{C}_R$ and $A_{(0,2)} \cap T (A_{(0,2)}) \cap T^2 (A_{(0,2)}) = \emptyset$. We continue this process for all $nm$ pairs in the index set $I_{(n,m)}$ (which corresponds to powers of $S$ and $T$) and obtain set $A_{(n-1,m-1)}$. Rename this set as $A_{(n,m)}$ to be consistent with statement of Theorem $\ref{Marker}$. The set $A_{(n,m)}$ satisfies $(i)$ and $(ii)$. This completes the proof. \hfill{$\square$}

With help of an example, we illustrate the proof of Theorem $\ref{Marker}$.

\begin{example} Assume $(n,m) = (2,3)$. We want to describe sets $A_{(i,j)}$'s ($0 \leq i < 2$, $0 \leq j < 3$), that we obtain at each step. The elements of $I_{(2,3)}$ in lexicographic order are $\{(0,0), (0,1),(0,2),(1,0),(1,1),(1,2)\}$ (they correspond to $S^0T^0 = e$, $T$, $T^2$, $S$, $ST$, $ST^2$ respectively). In the first step we obtain set $A_{(0,1)}$ such that $A_{(0,1)} \cap T A_{(0,1)} = \emptyset$. In the second step we obtain set $A_{(0,2)}$ such that the sets $A_{(0,2)}, T A_{(0,2)}, T^2 A_{(0,2)}$ are mutually disjoint.

The next element in the index set is $(1,0)$. So we have to go from $T^2$ to $S$, hence we would have to work with powers of $R$ that corresponds to $ST^{-2}$. In other words the partition $B_j$ will be $B_j = \{x \in  A_{(0,2)} : ST^{-1}x = R^j x\}$. This step will yield set $A_{(1,0)}$ such that sets 
\begin{equation*}
    \{A_{(1,0)}, ST^{-2} A_{(1,0)}, ST^{-1} A_{(1,0)}, S A_{(1,0)}\}
\end{equation*} are mutually disjoint. Similarly, the next step will yield set $A_{(1,1)}$ such that sets
\begin{equation*}
    \{A_{(1,1)}, T A_{(1,1)}, ST^{-1} A_{(1,1)}, S A_{(1,1)}, ST A_{(1,1)}\}
\end{equation*} are mutually disjoint. Finally we will obtain set $A_{(1,2)}$ such that sets
\begin{equation*}
\{A_{(1,2)}, T A_{(1,2)}, T^2 A_{(1,2)}, S A_{(1,2)}, ST A_{(1,2)}, ST^2 A_{(1,2)}\}
\end{equation*} are mutually disjoint. Denote $A_{(2,3)} = A_{(1,2)}$, thus sets $S^iT^j A_{(2,3)}$ are disjoint for $0 \leq i < 2$, $0 \leq j < 3$ as needed.
\end{example}

\noindent\textit{Proof of Theorem $\ref{d-Marker}$}: The proof of Theorem $\ref{d-Marker}$ is identical to the proof of Theorem $\ref{Marker}$. Instead of working with $2$-dimensional index set $I_{(n,m)}$, we will work with $d$-dimensional index set (with lexicographic ordering). Everything else remains the same.  \hfill{$\square$}

\section{Commuting automorphims and their coboundaries} \label{cob}

In this section, we provide proof for Theorem $\ref{thm main 1}$. As mentioned in the introduction our work is an extension of Kornfeld's ideas (see \cite[Theorem 1]{Kornfeld_1997}) to the domain of Borel dynamics. Although similar in nature our work differs from \cite{Kornfeld_1997} in the following manner. In \cite{Kornfeld_1997} the author used a version of $\Z^2$-Rokhlin Lemma (see \cite[Lemma 1]{Kornfeld_1997}) which he called a \textit{weak form} of Rokhlin Lemma for a measure preserving $\Z^2$-action (attributed to Conze \cite{Conze_72}). We use Theorem $\ref{Marker}$, which can be considered as a weak form of Rokhlin Lemma for free Borel $\Z^2$-actions. 

The other principal difference is that in \cite{Kornfeld_1997}, the author worked with an invariant ergodic probability measure on the ambient space. We do not work with any prescribed measure on the standard Borel space $(X,\B)$. However, in the proof of Theorem \ref{thm main 1} we use the existence of $\Gamma$-quasi-invariant probability measure on $(X,\B)$ (where $\Gamma$ is a free Borel $\Z^2$-action).

If $\Gamma$ is a Borel-action of $(X,\B)$ (and $E_X(\Gamma)$ denotes the corresponding orbit equivalence relation on $X$), we say that $\mu$ is $E_X(\Gamma)$-quasi invariant (or $\Gamma$-quasi invariant for short) if $\mu(B)> 0\, \Longleftrightarrow \mu(T(B)) > 0$ for all Borel sets $B \in \B$ and Borel maps $T: X \rightarrow X$ whose graphs are contained in $E_X(\Gamma)$. 

\begin{remark}\label{quasi} The following statement is a consequence of the fact that $\mu$ is $\Gamma$-quasi-invariant: Fix $\gamma \in \Gamma$ then for $\epsilon > 0$ there exists $\delta_{\gamma}$ such that for every $B \in \B$ with $\mu(B) < \delta_{\gamma}$ we have $\mu(\gamma B) < \epsilon$.

\end{remark}

\vspace{2mm}

\noindent\textit{Proof of Theorem $\ref{thm main 1}$.} By assumption $S$ and $T$ generate a free Borel $\Z^2$-action i.e. $T^mS^n x \neq x$, for all $x \in X$ and any pair of integers $(m,n) \neq (0,0)$. We will denote this Borel $\Z^2$-action by $\Gamma$ in rest of the proof. Also note that since $\Gamma$ is assumed to be free we can use Theorem $\ref{Marker}$. 

Without loss of generality we will construct a Borel function $f : X \rightarrow \mathbb{R}$, such that $f \in Cob(S)$ but $f \notin Cob(T)$. In particular, we will construct $f$ with following properties :

\noindent $(a)$ There exists a constant $M \in \mathbb{R}$, such that 

\begin{equation} \label{eq 5.1}
    \left |\underset{k=0}{\overset{n-1}{\sum}} f \circ S^{k}(x) \right | \leq M,
\end{equation} for every $x \in X$ and $n \in \mathbb{N}$. 

\noindent $(b)$ Let $\mu$ be a $\Gamma$-quasi-invariant probability measure on $(X,\B)$. For every $r\in \N$, there exists $m_r \in \mathbb{N}$ and a set $A_r \in \B$, $\mu(A_r) > \beta$ (for some $\beta >0$) such that 

\begin{equation}\label{eq 5.2}
    \left |\underset{k=0}{\overset{m_r - 1}{\sum}} f \circ T^{k} (x) \right | \geq r
\end{equation} for all $x \in A_r$.

Since $\Gamma$ is a free Borel $\Z^2$-action, in particular this implies that $T$ is not a power of $S$ hence $(\ref{eq 5.1})$ and $(\ref{eq 5.2})$ can hold simultaneously. Observe that $(\ref{eq 5.1})$ implies that $f \in Cob(S)$. To see this, set 

\begin{equation*}
    g(x) = \underset{n \geq 1}{\textrm{sup}} \, \Bigg(\underset{k=0}{\overset{n-1}{\sum}} f \circ S^{k} (x)\Bigg).
\end{equation*}  Thus $f(x) = g(x)- g(Sx)$ for all $x \in X$.

We claim that property $(b)$ implies that $f \notin Cob(T)$. To see this assume by contradiction that $f \in Cob(T)$. Thus, there exists a transfer function $g(x)$ such that $f(x) = g(x)- g(Tx)$. For $n\in \N$, we can write $\underset{k=0}{\overset{n-1}{\sum}} f \circ T^{k}(x) = g(x) - g \circ T^{n} (x)$.  We will use the statement in Remark $\ref{quasi}$ as follows : Fix $\gamma \in \Gamma$ then for $\beta > 0$ (where $\beta$ is as in $(b)$ above) and every $0 <\epsilon < \beta$, there exists a $\delta_{\gamma} > 0$, such that $\epsilon + \delta_{\gamma} < \beta$ and for every $B \in \B$ with $\mu(B) < \delta_{\gamma}$, we have $\mu(\gamma B) < \epsilon$. 

For $\gamma = T^{-1}$ and $\epsilon >0$ such that $0 <\epsilon < \beta$ (where $\beta$ is as in $(b)$ above) and consider $\delta_{T^{-1}}$ (write $\delta_1$ for simplicity) as described above. For $K_1 \in \R$ such that 
\begin{equation}\label{eq K}
    \mu(\{x : |g(x)| \geq K_1 \}) < \delta_1,
\end{equation} using the fact that $\mu$ is $\Gamma$-quasi-invariant we get $\mu(\{x : |g(T x)| \geq K_1 \}) < \epsilon$. We choose $\delta_1$ to be such that $0< \epsilon + \delta_1 < \beta$ and $K_1$ to be the infimum possible corresponding to $\delta_1$. Thus we get 
\begin{equation*}
\mu\Big( \Big\{x : \left |g(x) - g(Tx) \right | \geq 2K_1  \Big\}  \Big) < \epsilon + \delta_1 < \beta.
\end{equation*} i.e. $\mu( \{x : \left |f(x)\right | \geq 2K_1\})< \beta$. Now we repeat the process for every $n > 1$ (put $\gamma = T^{-n}$).  For same $\epsilon >0 $ as above we obtain a sequence of $\delta_n >0$. For every $n > 1$, we choose each $\delta_n$ such that $0<\epsilon + \delta_n < \beta$ and $\delta_{n} < \delta_{n-1}$. Thus for $K_n \in \R$ such that 
\begin{equation*}
    \mu(\{x : |g(x)| \geq K_n \}) < \delta_n,
\end{equation*} we get, \begin{equation*}
    \mu(\{x : |g(T^n x)| \geq K_n \}) < \epsilon
\end{equation*} for every $n>1$. As before we work with infimum possible $K_n$ corresponding to $\delta_n$. Since $\{\delta_n\}_{n=1}^{\infty}$ is a decreasing sequence we get that the corresponding sequence $\{K_n\}_{n=1}^{\infty}$ is increasing. Hence for every $n>1$ we obtain 
\begin{equation*}
  \mu\Big( \Big\{x : \left |g(x) - g(T^{n} x) \right | \geq 2K_n  \Big\}  \Big) < \epsilon + \delta_n < \beta.  
\end{equation*} Thus for every $n \in \N$ we obtain $\delta_n$ and $K_n$ such that,
\begin{equation}
    \mu\Big( \Big\{x : \left |\underset{k=0}{\overset{n-1}{\sum}} f \circ T^{k} (x)\right | \geq 2K_n  \Big\}  \Big) < \epsilon + \delta_n < \beta.
\end{equation} Assume that $\underset{n}{\mathrm{sup}}\,\, K_n = : K$ exists (we will prove this later). Thus for every $n \in \N$ \begin{equation*}
    \mu(\{x : |g(x)| \geq K \}) < \delta_n,
\end{equation*} we get \begin{equation*}
    \mu(\{x : |g(T^n x)| \geq K \}) < \epsilon.
\end{equation*} Thus for every $n \in \N$ 
\begin{equation}
    \mu\Big( \Big\{x : \left |\underset{k=0}{\overset{n-1}{\sum}} f \circ T^{k} (x)\right | \geq 2K  \Big\}  \Big) < \epsilon + \delta_n < \beta.
\end{equation} which contradicts $(\ref{eq 5.2})$. We need to show that $K : = \underset{n}{\mathrm{sup}}\,\, K_n$ exists. To see this assume by contradiction that the supremum does not exists. Hence for every $M \in \R_{+}$ there exists $n \in \N$ such that, \begin{equation*}
    \mu(\{x : |g(x)| \geq M/2 \}) > \delta_n.
\end{equation*} Since $\mu$ is $\Gamma$-quasi-invariant we get
\begin{equation*}
    \mu(\{x : |g(T x)| \geq M/2 \}) > 0.
\end{equation*} Hence $ \mu(\{x : |g(x) - g(T x)| \geq M \}) > 0$, i.e. $ \mu(\{x : |f(x)| \geq M \}) > 0$. Since we can do this for every $M \in \R_{+}$, $f$ is essentially unbounded, which is a contradiction.

Now we construct a Borel function $f$ with properties $(a)$, $(b)$.   The function $f$ will be constructed as the sum of infinite series $f := \underset{r=1}{\overset{\infty}{\sum}} f_r $, in which every term $f_r$ is associated with a certain Rokhlin tower $\xi_r$ of the $\mathbb{Z}^2$-action generated by commuting Borel automorphisms $S$ and $T$. Below we describe the construction of tower $\xi_r$ and the associated function $f_r$.

The size of tower $\xi_r$ (given by $(n_rm_r)$) is determined by two increasing sequence $\{n_r\}$ and $\{m_r\}$ of natural numbers. We also associate a decreasing sequence of positive real numbers $\{\alpha_r\}$ with towers $\xi_r$ (we assume $\alpha_r \rightarrow
 0$). The only restriction on sequence $\{n_r\}$ is that it is an increasing sequence of natural numbers. On the other hand, we have the following assumptions for sequence $\{m_r\}$ and $\{\alpha_r\}$:
\begin{equation}\label{eq 5.3}
    m_r \Big(\underset{s=r+1}{\overset{\infty}{\sum}} \alpha_s \Big) \leq 1, \hspace{15mm} r = \{1,2,...\}
\end{equation}

\begin{equation}\label{eq 5.4}
    m_r \alpha_r \geq 2r + \Big(\underset{t=1}{\overset{r-1}{\sum}} \alpha_t m_t \Big), \,\,\,\, r = \{2,3,...\}
\end{equation} 
\begin{figure}[h]
\centering
\includegraphics[width=60mm]{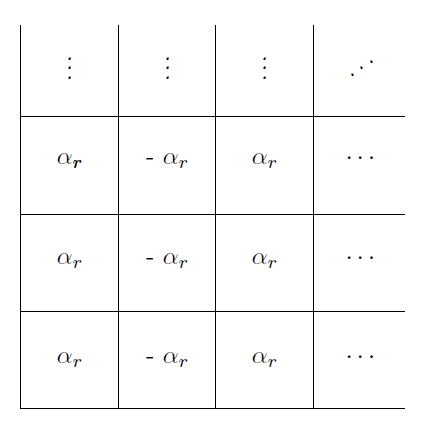}
\caption{The $\Z^2$ tower $\xi_r$.}
\label{fig:tower}
\end{figure}

The sequences $\{m_r\}$ and $\{\alpha_r\}$ with above properties can be constructed inductively. We use Theorem $\ref{Marker}$ to define tower $\xi_r$ as follows. Set $n = n_r$ and $m = m_r$ in Theorem $\ref{Marker}$ to obtain a set $A_r = A_{(n,m)}$, which is the base of tower $\xi_r$ (the bottom left block in Figure \ref{fig:tower}). This tower is a rectangle made up of $n_r m_r$ disjoint blocks each representing a set of the form $S^iT^j A_r$, where $0\leq i < n_r$ and $0\leq j < m_r$ (see Figure \ref{fig:tower}). The horizontal direction in the tower corresponds to the transformation $S$ and the vertical corresponds to the transformation $T$.

 Function $f_r$ is zero outside the tower $\xi_r$. On the tower $\xi_r$ it is defined to be constant on each square (in other words $f_r$ is constant on each set $S^iT^j A_r \subset \xi_r$). The value of $f_r$ on each square is defined as follows : In each row the value alternates $ \alpha_r$, $- \alpha_r$ starting with plus sign on the leftmost set. (see Figure \ref{fig:tower}).  
 
We now estimate the sums $\left |\underset{k=0}{\overset{n-1}{\sum}} f_r \circ S^{k} (x) \right |$ and $\left |\underset{k=0}{\overset{n-1}{\sum}} f_r \circ T^{k} (x) \right |$ for fixed $r$. The first sum will be estimated from below to prove that $f$ satisfies property $(a)$ and the second sum will be estimated from above (for $x \in A_r$) to show that $f$ satisfies property $(b)$. Note for any $x \in \xi_r$ we have  $\left |\underset{k=0}{\overset{n-1}{\sum}} f_r \circ S^{k} (x) \right | \leq \alpha_r$. By definition $f_r = 0$ outside the tower $\xi_r$, thus we get
\begin{equation*}
    \left |\underset{k=0}{\overset{n-1}{\sum}} f_r \circ S^{k} (x) \right | \leq \alpha_r \,\,\, \textrm{for}\,\,\, x \in X.
\end{equation*} Hence for $x \in X$,
\begin{equation*}
     \left |\underset{k=0}{\overset{n-1}{\sum}} f \circ S^{k} (x) \right | =  \left |\underset{k=0}{\overset{n-1}{\sum}} \Bigg(\underset{r=1}{\overset{\infty}{\sum}} f_r\Bigg) \circ S^{k} (x) \right | \leq \underset{r=1}{\overset{\infty}{\sum}}\left |\underset{k=0}{\overset{n-1}{\sum}} f_r \circ S^{k} (x) \right | \leq \underset{r=1}{\overset{\infty}{\sum}}  \alpha_r =: M
\end{equation*} Therefore $f$ satisfies property $(a)$. 

To see that $f$ satisfies $(b)$, note that to every $r \in \mathbb{N}$ we have associated a pair of natural number $(n_r, m_r)$ and a tower $\xi_r$ with base $A_r$. Since $A_r$ is a complete section with respect to $\Gamma$ and $\mu$ is a $\Gamma$-quasi-invariant probability measure, $\mu(A_r)> 0$. We work with $A_r$ and the corresponding tower $\xi_r$. Note that $\left |\underset{k=0}{\overset{m_r - 1}{\sum}} f_r \circ T^{k} (x)\right| = m_r \alpha_r$, for all $x \in A_r$. Also we can write
\begin{equation}\label{eq 5.5}
    \left |\underset{k=0}{\overset{m_r - 1}{\sum}} f \circ T^{k} (x)\right | \geq \left |\underset{k=0}{\overset{m_r - 1}{\sum}} f_r \circ T^{k} (x) \right | - \underset{t=1}{\overset{r-1}{\sum}} \left |\underset{k=0}{\overset{m_r - 1}{\sum}} f_t \circ T^{k} (x) \right | - \underset{s=r+1}{\overset{\infty}{\sum}} \left |\underset{k=0}{\overset{m_r - 1}{\sum}} f_s \circ T^{k} (x) \right |.
\end{equation} Since $f_t$ is zero outside the tower $\xi_t$ and $t < r$ we get
\begin{equation}\label{t}
\left |\underset{k=0}{\overset{m_r - 1}{\sum}} f_t \circ T^{k} (x) \right | = \left |\underset{k=0}{\overset{m_t - 1}{\sum}} f_t \circ T^{k} (x) \right | \leq  \alpha_t m_t.
\end{equation} for $x \in A_r$. In $(\ref{t})$ if towers $\xi_r$ and $\xi_t$ are disjoint then the sum is zero. If towers $\xi_r$ and $\xi_t$ are not disjoint than the sum can be at most $\alpha_t m_t$. This implies for $x \in A_r$,
\begin{equation*}
- \underset{t=1}{\overset{r-1}{\sum}} \left |\underset{k=0}{\overset{m_r - 1}{\sum}} f_t \circ T^{k} (x) \right | \geq - \underset{t=1}{\overset{r-1}{\sum}} \alpha_t m_t .
\end{equation*} Since $f_s$ is zero outside the tower $\xi_s$ and $s > r$ similarly we get
\begin{equation*}
\left |\underset{k=0}{\overset{m_r - 1}{\sum}} f_s \circ T^{k} (x) \right |  \leq  \alpha_s m_r
\end{equation*} for $x \in A_r$. This implies 
\begin{equation*}
    - \underset{s=r+1}{\overset{\infty}{\sum}} \left |\underset{k=0}{\overset{m_r - 1}{\sum}} f_s \circ T^{k} (x) \right | \geq - \underset{s=r+1}{\overset{\infty}{\sum}} \alpha_s m_r
\end{equation*} for $x \in A_r$. Thus by $(\ref{eq 5.5})$ 
\begin{equation}\label{eq 5.6}
\left |\underset{k=0}{\overset{m_r - 1}{\sum}} f \circ T^{k} (x) \right | \geq m_r \alpha_r -  \underset{t=1}{\overset{r-1}{\sum}} \alpha_t m_t -  \underset{s=r+1}{\overset{\infty}{\sum}} \alpha_s m_r
\end{equation} for $x \in A_r$. Hence by $(\ref{eq 5.4})$ and $(\ref{eq 5.3})$ 
\begin{equation}\label{eq 5.7}
 \left |\underset{k=0}{\overset{m_r - 1}{\sum}} f \circ T^{k} (x) \right | \geq 2r -  \underset{s=r+1}{\overset{\infty}{\sum}} \alpha_s m_r \geq 2r - 1 \geq r
\end{equation} for $x \in A_r$. This shows that $f$ satisfies property $(b)$, which completes the proof. 

\hfill{$\square$}
 
\textbf{Acknowledgments.} The author is very thankful to Sergii Bezuglyi and Palle E. T. Jørgensen for guidance during this work. The author would also like to thank the referee for careful reading and valuable suggestions.

\bibliographystyle{plain}
\bibliography{references1.bib}

\noindent Shrey Sanadhya\\  
shrey-sanadhya@uiowa.edu\\

\noindent {\small
\noindent The University of Iowa\\
Department of Mathematics\\
14 MacLean Hall,
Iowa City, Iowa 52242, USA
}\bigskip

\end{document}